\documentclass[12pt]{article}
\usepackage{amsmath}
\usepackage{amssymb}
\allowdisplaybreaks

\def\q{\quad}
\def\qq{\qquad}

\def\mod#1{\ (\text{\rm mod}\ #1)}
\def\t{\text}
\def\f{\frac}
\def\e{\equiv}
\def\b{\binom}

 \def\ls#1#2{\big(\f{#1}{#2}\big)}
\def\Ls#1#2{\Big(\f{#1}{#2}\Big)}

\def\s3{(-1)^{[\frac p3]}}

\begin{document}
 \centerline {\bf
Identities and congruences involving orthogonal polynomials and Ap\'ery-like numbers}
\par\q\newline
\centerline{Zhi-Hong Sun}\newline \centerline{School of Mathematics
and Statistics} \centerline{Huaiyin Normal University}
\centerline{Huaian, Jiangsu 223300, P.R. China} \centerline{Email:
zhsun@hytc.edu.cn} \centerline{URL:
https://maths.hytc.edu.cn/szh.htm}
\vskip0.2cm
\par\q
\par {\bf Abstract.} In this paper, we establish a general identity for three-term recurrence sequences and then give applications to orthogonal polynomials and Ap\'ery-like numbers.
\par\q

\section*{1. Introduction}

\par\q Many orthogonal polynomials and Ap\'ery-like numbers $\{u_n\}$ satisfy $u_0=1$ and the three-term recurrence relation:
$$(n+1)^ru_{n+1}=P(n)u_n-cn^ru_{n-1}\q(n=0,1,2,\ldots),$$
where $c\not=0$, $r\in\{0,1,2,3\}$ and $P(n)$ is a polynomial of $n$ with degree less than $4$. For instance, the Legendre polynomials $\{P_n(x)\}$, Laguerre polynomials $\{L_n(x)\}$,  Bessel polynomials $\{y_n(x)\}$ and the Bessel function $J_{n}(x)$ satisfy the following relations:
\begin{align*}
&P_{-1}(x)=0,\ P_0(x)=1,\ (n+1)P_{n+1}(x)=(2n+1)xP_n(x)-nP_{n-1}(x)\ (n\ge 0),
\\&L_{-1}(x)=0,\ L_0(x)=1,\ (n+1)L_{n+1}(x)=(2n+1-x)L_n(x)-nL_{n-1}(x) \ (n\ge 0),
\\&y_{-1}(x)=1,\ y_0(x)=1,\ y_{n+1}(x)=(2n+1)xy_n(x)+y_{n-1}(x)\ (n\ge 0),
\\&J_{n+1}(x)=\f{2n}xJ_n(x)-J_{n-1}(x)\ (n\ge 0).
\end{align*}
For the properties of these special functions, see [1,4,11,27].
\par Let $\Bbb Z,\Bbb Z^+$ and $\Bbb R$ be the set of integers, the set of positive integers and the set of real numbers, respectively. For $a,b,c\in\Bbb R$ with $c\not=0$ let $\{u_n\}$ be given by
\begin{align*}u_0=1, \ u_1=\ b,\ (n+1)^3u_{n+1}
=(2n+1)(an(n+1)+b)u_n-cn^3u_{n-1}\  (n\ge 1).\tag {1.1}\end{align*}
 If $a,b,c\in\Bbb Z$ and $u_n\in\Bbb Z$ for $n=1,2,3,\ldots$, following [20] and [22] we call $\{u_n\}$ the first kind of Ap\'ery-like numbers. Let
\begin{align*}
 &A_n= \sum_{k=0}^n\binom nk^2\binom{n+k}k^2,
\q D_n=\sum_{k=0}^n\b nk^2\b{2k}k\b{2n-2k}{n-k},\\&
b_n=\sum_{k=0}^{[n/3]}\b{2k}k\b{3k}k\b n{3k}\b{n+k}k(-3)^{n-3k}, \
T_n=\sum_{k=0}^n\b nk^2\b{2k}n^2, \
\\\qq\qq& V_n=\sum_{k=0}^n\b
nk\b{n+k}k(-1)^k\b{2k}k^216^{n-k}
 =\sum_{k=0}^n\b{2k}k^2\b{2n-2k}{n-k}^2,
 \\&V_n^{(3)}=\sum_{k=0}^n\b
nk\b{n+k}k(-1)^k\b{2k}k\b{3k}k27^{n-k},
\\&V_n^{(4)}=\sum_{k=0}^n\b
nk\b{n+k}k(-1)^k\b{2k}k\b{4k}{2k}64^{n-k},
\\&V_n^{(6)}=\sum_{k=0}^n\b
nk\b{n+k}k(-1)^k\b{3k}k\b{6k}{3k}432^{n-k} .
\end{align*}
Then $A_n,D_n,b_n,T_n,V_n,V_n^{(3)},V_n^{(4)}$ and $V_n^{(6)}$ are the first kind of Ap\'ery-like numbers with
\begin{align*}(a,b,c)=&(17,5,1),(10,4,64),(-7,-3,81),
(12,4,16),\\&(16,8,256),
(27,15,729),(64,40,4096),(432,312,186624),\end{align*}
respectively. The numbers $\{A_n\}$, $\{D_n\}$ and $\{b_n\}$ are called Ap\'ery numbers, Domb numbers and Almkvist-Zudilin numbers, respectively. For the discussions on $\{A_n\}$, $\{D_n\}$, $\{b_n\}$, $\{T_n\}$ and
$\{V_n\}$ see [20], [22] and A005259, A002895, A125143, A290575 and A036917 in [18].
\par For $a,b,c\in\Bbb R$ with $c\not=0$ let $\{u_n\}$ be given by
\begin{align*}u_0=1, \ u_1=\ b,\ (n+1)^2u_{n+1}
=(an(n+1)+b)u_n-cn^2u_{n-1}\  (n\ge 1).\tag{1.2}\end{align*}
 If $a,b,c\in\Bbb Z$ and $u_n\in\Bbb Z$ for $n=1,2,3,\ldots$, we call $\{u_n\}$ the second kind of Ap\'ery-like numbers. Let
 \begin{align*}
 &A'_n=\sum_{k=0}^n\b nk^2\b {n+k}k,
 \q f_n=\sum_{k=0}^n\b nk^3\q
 \\&S_n=\sum_{k=0}^n\b nk\b{2k}k\b{2n-2k}{n-k},
 \q a_n=\sum_{k=0}^n\binom nk^2\binom{2k}k,
  \\&W_n=\sum_{k=0}^{[n/3]}\binom{2k}k\binom{3k}k\binom n{3k}(-3)^{n-3k},
 \q Q_n=\sum_{k=0}^n\binom nk(-8)^{n-k}f_k\end{align*}
and
  \begin{align*}&G_n=\sum_{k=0}^n\binom
nk(-1)^k\binom{2k}k^216^{n-k},
\\&G_n^{(3)}=\sum_{k=0}^n\binom
nk(-1)^k\binom{2k}k\binom{3k}k27^{n-k},
\\&G_n^{(4)}=\sum_{k=0}^n\binom
nk(-1)^k\binom{2k}k\b{4k}{2k}64^{n-k},
\\&G_n^{(6)}=\sum_{k=0}^n\binom
nk(-1)^k\binom{3k}k\binom{6k}{3k}432^{n-k},\end{align*}
where $[x]$ is the greatest integer not exceeding $x$.
Then $A'_n, f_n, S_n, a_n, W_n,
 Q_n, G_n,G_n^{(3)}, G_n^{(4)}$ and $G_n^{(6)}$
 are the second kind of Ap\'ery-like sequences with
\begin{align*}(a,b,c)=&(11,3,-1),(7,2,-8),(12,4,32),(10,3,9),
(-9,-3,27), (-17,-6,72),
\\&(32,12,256),(54,21,729),(128,52,4096),(864,372,186624),
\end{align*} respectively. For the properties of $G_n,G_n^{(3)}, G_n^{(4)}$ and $G_n^{(6)}$ see [21]. The  sequence
 $\{f_n\}$ is called Franel numbers.  For the discussions on $\{A'_n\},\ \{f_n\},$ $\{S_n\},\ \{a_n\},\
 \{Q_n\}$
and $\{W_n\}$, see A005258, A000172, A081085, A002893, A093388 and A291898 in Sloane's database [18].

\par Recently, Mao and Ni [15] established useful identities for Ap\'ery-like sequences. If $\{u_n\}$ is the sequence given by (1.2), they showed that for $x\not=0$ and $p\in\Bbb Z^+$,
$$\sum_{n=0}^{p-1}\big((x^2-ax+c)n^2+(2c-ax)n+c-bx\big)\f{u_n}{x^n}
=\f{p^2}{x^{p-1}}(cu_{p-1}-xu_p).\eqno{(1.3)}$$
If $\{u_n\}$ is the sequence satisfying (1.1), they proved that for $x\not=0$ and $p\in\Bbb Z^+$,
 \begin{align*}&\sum_{n=0}^{p-1}\big((x^2-2ax+c)n^3+(3c-3ax)n^2+(3c-(a+2b)x)n+c-bx\big)\f{u_n}{x^n}
\\&=\f{p^3}{x^{p-1}}(cu_{p-1}-xu_p).\tag{1.4}\end{align*}
\par For the Legendre polynomials $\{P_n(x)\}$, it is well known that $(1-x^2)P'_n(x)=nP_{n-1}(x)-nxP_n(x)\ (n\ge 1)$.
 Thus, from the Christoffel-Darboux formula one may deduce the formula (see [2])
$$\sum_{n=0}^{p-1}(2n+1)P_n(x)^2
=\f{p^2}{1-x^2}\big(P_p(x)^2-2xP_{p-1}(x)P_p(x)+P_{p-1}(x)^2\big).
\eqno{(1.5)}$$
 \par In [24] and [5], Sun and Li investigated the properties of Ap\'ery-like polynomials $g_n(x)$ and $v_n(x)$ defined by
\begin{align*}
&g_0(x)=1,\ g_1(x)=\frac{x+1}2,\\& (n+1)^2g_{n+1}(x)=\Big(2n(n+1)+\frac{x+1}2\Big)g_n(x)-n^2g_{n-1}(x)\ (n\ge 1),
\\&v_0(x)=1,\ v_1(x)=x,\\& (n+1)^3v_{n+1}(x)=(2n+1)(n(n+1)+x)v_n(x)-n^3v_{n-1}(x)\ (n\ge 1).\end{align*}
In particular, Li and Sun [5] discovered the identities
\begin{align*}
&\sum_{n=0}^{p-1}(2n+1)g_n(x)^2
= p^2g_p(x)g_{p-1}(x)-\frac{2p^4}{x-1}\bigl(g_p(x)
-g_{p-1}(x)\bigr)^2\tag{1.6}
\end{align*}
and
\begin{align*}\sum_{n=0}^{p-1}\big((2n+1)^3+(2n+1)x\big)v_n(x)^2
 =-\f{p^6}{x-1}(v_p(x)-v_{p-1}(x))^2+2p^4v_p(x)v_{p-1}(x).
 \tag{1.7} \end{align*}
\par Motivated by (1.3)-(1.7), in this paper we construct a very general identity for
the sequence $\{u_n\}$ satisfying
$$(n+1)^ru_{n+1}=P(n)u_n-cn^ru_{n-1}\q(n=0,1,2,\ldots),
\eqno{(1.8)}$$
where $c\not=0$, $r\ge 0$ and $P(n)$ is a polynomial of $n$ with degree less than $4$. See Theorem 2.1. As consequences and applications, we deduce many general identities for  orthogonal polynomials and Ap\'ery-like numbers, and deduce many supercongruences involving Ap\'ery-like numbers. Here are some typical results:
\begin{align*}
&\sum_{n=0}^{p-1}\big(9(2n+1)^5-4(2n+1)^3-
(2n+1)\big)A_n^2\e -20p^4+24p^6+\f{16}3p^7B_{p-3}\mod {p^8},
\\&\sum_{n=0}^{p-1}\big(3(2n+1)^5-2(2n+1)^3
-(2n+1)\big)\frac{T_n^2}{16^n}
\e -8p^4+8p^6\mod {p^8},
\\&\sum_{n=0}^{p-1}\big(25(2n+1)^3+3(2n+1)\big)(-1)^n{A'_n}^2
\e -12p^2+40p^4-\f{40}3p^5B_{p-3}\mod {p^6},
\\&\sum_{n=0}^{p-1}(18n^3+75n^2+78n+25)
\frac{S_n^2}{16^n}
\\&\e  5p^4-4(-1)^{\f{p-1}2}p^2(9p^2-3p-1)2^{p-1}+12p^4E_{p-3}
\mod {p^5},
\\&\sum_{n=0}^{p-1}(2646n^3+5103n^2+3668n+940)f_n^2
\\&\e 32p^2(21p-10)8^{p-1}+196p^4(3+8(8^{p-1}-1))
-280p^5B_{p-3} \mod {p^6},
\end{align*}
where $p$ is a prime greater than $3$, $\{B_n\}$ and $\{E_n\}$ are the Bernoulli numbers and Euler numbers given by
\begin{align*} &B_0=1,\q\sum_{k=0}^{n-1}\b nkB_k=0\q(n\ge 2),\\& E_{2n-1}=0,\q E_0=1\q\t{and}\q E_{2n}=-\sum_{k=0}^{n-1}
\b {2n}{2k}E_{2k}\q(n\ge 1).
\end{align*}

\section*{2. A general identity}
\par\q In this section, we establish a general identity, which can be viewed as a vast generalization of (1.5)-(1.7).
\vskip0.2cm
\par{\bf Theorem 2.1} {\sl Let $c,x\in\Bbb R$ with $cx(x+1)\not=0$, and $P (n)=a_3n^3+a_2n^2+a_1n+a_0$. Suppose that $r\ge 0$ and $\{u_n\}$ satisfies the recurrence relation:
$$(n+1)^ru_{n+1}=P(n)u_n-cn^ru_{n-1}\ (n=0,1,2,\ldots).$$
Then
\begin{align*}&\sum_{n=0}^{p-1}P_r(n,x)\f{u_n^2}{(cx)^n}
-\f 1{(cx)^{p-1}}\Big(p^{2r}u_p^2+\f cxp^{2r}u_{p-1}^2-p^rs_p(x)u_pu_{p-1}\Big)
\\&=\begin{cases}
0&\t{if $r>0$,}
\\-cx(u_0^2+\f cxu_{-1}^2-s_0(x)u_0u_{-1})&\t{if $r=0$,}
\end{cases}\end{align*}
where
\begin{align*}
s_n(x) = & \frac{2a_3}{x+1} n^3 + \Big(\frac{2a_2}{x+1} -\frac{6a_3} {(x+1)^2}\Big)n^2 +\Big ( \frac{2a_1}{x+1} - \frac{4a_2}{(x+1)^2} -
\frac{6a_3 (x-1)}{(x+1)^3}\Big) n \\
&+ \frac{2a_0}{x+1} - \frac{2a_1}{(x+1)^2} -\frac{2a_2(x-1)}{(x+1)^3}
-\frac{2a_3(x^2-4x+1)}{(x+1)^4}
\end{align*}
and
$$P_r(n,x)=\frac cx(n+1)^{2r}-cxn^{2r}-P(n)^2+xs_{n}(x)P(n).$$}
\par{\it Proof.} It is easy to check the identity
$$s_{n+1}(x)+xs_n(x)=2P(n)\ (n\ge 0).$$
Set
$$h_r(n)=\f 1{(cx)^{n-1}}\Big(n^{2r}u_n^2+\f cxn^{2r}u_{n-1}^2 -n^rs_n(x)u_nu_{n-1}\Big)\ (n\ge 0).$$
Then
\begin{align*}
&h_r(n+1)\\&=\f 1{(cx)^{n}}\big((n+1)^{2r}u_{n+1}^2+\f cx(n+1)^{2r}u_n^2
-s_{n+1}(x)(n+1)^ru_{n+1}u_n\big)
\\&=\f 1{(cx)^{n}}
\big((P(n)u_n-cn^ru_{n-1})^2+\f cx(n+1)^{2r}u_n^2
-s_{n+1}(x)
(P(n)u_n-cn^ru_{n-1})u_n\big)
\\&=\f 1{(cx)^{n}}
\Big(\big(P(n)^2+\f cx(n+1)^{2r}-s_{n+1}(x)P(n)\big)u_n^2
\\&\q+c^2n^{2r}u_{n-1}^2+\big(-2cn^rP(n)
+cn^rs_{n+1}(x)\big)u_nu_{n-1}
\Big)
\\&=h_r(n)+\f 1{(cx)^{n}}
\Big(\big(P(n)^2+\f cx(n+1)^{2r}-s_{n+1}(x)P(n)-cxn^{2r}
\big)u_n^2
\\&\q+cn^r\big(-2P(n)
+s_{n+1}(x)+xs_n(x)\big)u_nu_{n-1}
\Big)
\\&=h_r(n)+\f 1{(cx)^{n}}
\Big(\big(P(n)^2+\f cx(n+1)^{2r}-2P(n)^2+xs_n(x)P(n)-cxn^{2r}
\big)u_n^2
\\&=h_r(n)+P_r(n,x)\f{u_n^2}{(cx)^{n}}.
\end{align*}
Thus,
\begin{align*}\sum_{n=0}^{p-1}P_r(n,x)\f{u_n^2}{(cx)^n}
&=\sum_{n=0}^{p-1}(h_r(n+1)-h_r(n))=h_r(p)-h_r(0)
\\&=\begin{cases}h_r(p)&\t{if $r>0$,}
\\h_r(p)-h_r(0)&\t{if $r=0$.}\end{cases}
\end{align*}
This proves the theorem.
 \vskip0.2cm
 \par{\bf Remark 2.1} By the proof of Theorem 2.1, replacing $n^r$ and $(n+1)^r$ in Theorem 2.1 with arbitrary functions $t(n)$ and $t(n+1)$ one can obtain a more general result.
 \vskip0.2cm
 \par As consequences, from Theorem 2.1 we have the following formulas for the Legendre polynomials $P_n(x)$, Bessel polynomials $y_n(x)$, Laguerre polynomials $L_n(x)$, the polynomial $d_n(x)=\sum_{k=0}^n\b nk\b xk2^k$ and the Bessel function $J_n(x)$.
 \vskip0.2cm
 \par{\bf Corollary 2.1} {\sl Let $\{P_n(t)\}$ be the Legendre
 polynomials, $x\in\Bbb R$ and $x\not=0,-1$. Then
 \begin{align*}&\sum_{n=0}^{p-1}\Big(\Big(\frac{1-x^2}x
 +\frac{4(x-1)}{x+1}t^2\Big)n^2
 +\Big(\frac 2x+\frac{4(x^2-2x-1)t^2}{(x+1)^2}\Big)n+
 \frac 1x+\frac{(x^2-4x-1)t^2}{(x+1)^2}\Big)\frac{P_n(t)^2}{x^n}
 \\&=\frac 1{x^{p-1}}\Big(p^2P_p(t)^2+\frac{p^2}xP_{p-1}(t)^2
 -\frac{2pt(2p(x+1)+x-1)}{(x+1)^2}P_p(t)P_{p-1}(t)\Big).
 \end{align*}
 Hence, for the Delannoy numbers $D_n=\sum_{k=0}^n\b nk\b {n+k}k=P_n(3)$ we have
 \begin{align*}
&\sum_{n=0}^{p-1}\big((1-x^2)(x^2-34x+1)n^2
+(36x^3-70x^2-32x+2)n
+9x^3-35x^2-7x+1\big)\frac{D_n^2}{x^n}\\
&\quad =\frac{1}{x^{p-1}}\big(
x(x+1)^2p^2D_p^2+(x+1)^2p^2D_{p-1}^2
-6px\big((2p+1)x+2p-1\big)D_pD_{p-1}\big),
\end{align*}
 for $b,c\in\Bbb R$ with $b^2-4c\not=0$ and
 $$T_n(b,c)=\sum_{k=0}^{[n/2]}\binom n{2k}\binom{2k}kb^{n-2k}c^k
 =(\sqrt{b^2-4c})^nP_n\Ls b{\sqrt{b^2-4c}}$$
 we have
  \begin{align*}&\sum_{n=0}^{p-1}\Big(\Big(\frac{1-x^2}x
 +\frac{4(x-1)}{x+1}\cdot
 \frac{b^2}{b^2-4c}\Big)n^2
 +\Big(\frac 2x+\frac{4(x^2-2x-1)}{(x+1)^2}\cdot\frac {b^2}{b^2-4c}\Big)n
 \\&\quad+ \frac 1x+\frac{(x^2-4x-1)}{(x+1)^2}\cdot
 \frac{b^2}{b^2-4c}\Big)\frac{T_n(b,c)^2}
 {((b^2-4c)x)^n}
 \\&=\frac p{((b^2-4c)x)^{p-1}}\Big(p\frac{T_p(b,c)^2}{b^2-4c}+
 \frac{p}xT_{p-1}(b,c)^2
 -\frac{2b(2p(x+1)+x-1)}{(b^2-4c)(x+1)^2}
 T_p(b,c)T_{p-1}(b,c)\Big).
 \end{align*}
 }
 \par{\it Proof.} Taking $c=1$, $r=1$, $P(n)=(2n+1)t$, $a_0=t,\ a_1=2t,\ a_2=a_3=0$, $u_n=P_n(t)$ in Theorem 2.1 gives
 $$s_n(x)=\f{4tn}{x+1}+\f{2(x-1)t}{(x+1)^2}$$
 and
 \begin{align*}P_1(n,x)&=\f{(n+1)^2}x-n^2x-(2n+1)^2t^2+\f{2tx}{x+1}\Big(2n+\f{x-1}{x+1}
 \Big)(2n+1)t
 \\&=\Big(\f{1-x^2}x+\f{4(x-1)}{x+1}t^2\Big)n^2
 +\Big(\f 2x+\f{4(x^2-2x-1)t^2}{(x+1)^2}\Big)n+
 \f 1x+\f{(x^2-4x-1)t^2}{(x+1)^2}.\end{align*}
 Thus, the result for $P_n(t)$ follows from Theorem 2.1.
 Taking $t=3$ yields the formula for Delannoy numbers $D_n$, and taking $t=\f b{\sqrt{b^2-4c}}$ yields the formula for $T_n(b,c)$.
\vskip0.2cm
\par{\bf Corollary 2.2} {\sl Let $\{y_n(t)\}$ be the Bessel
 polynomials and $x\in\Bbb R$ with $x\not=0,-1$. Then}
 \begin{align*}&\sum_{n=0}^{p-1}\Big(\frac{4(x-1)t^2}{x+1}n^2
 +\frac{4(x^2-2x-1)t^2}{(x+1)^2}n
 +\frac{x^2-1}x+\frac{(x^2-4x-1)t^2}{(x+1)^2}\Big)\frac
 {y_n(t)^2}{(-x)^n}
 \\&=\frac 1{(-x)^{p-1}}
 \Big(y_p(t)^2-\f 1xy_{p-1}(t)^2
-\frac{2(2p(x+1)+x-1)t}{(x+1)^2}y_p(t)y_{p-1}(t)\Big)
+x-1-\frac{2x(x-1)t}{(x+1)^2}.
 \end{align*}
 \par{\it Proof.} Taking $c=-1$, $r=0$, $P(n)=(2n+1)t$, $a_0=t,\ a_1=2t,\ a_2=a_3=0$, $u_n=y_n(t)$ in Theorem 2.1 gives
 $$s_n(x)=\f{4tn}{x+1}+\f{2(x-1)t}{(x+1)^2}
 =\f{2t(2n(x+1)+x-1)
 }{(x+1)^2}$$
 and
 \begin{align*}P_0(n,x)&=-\f 1x+x-(2n+1)^2t^2+\f{2tx(2n(x+1)+x-1)
 }{(x+1)^2}(2n+1)t
 \\&=\f{4(x-1)t^2}{x+1}n^2+\f{4(x^2-2x-1)t^2}{(x+1)^2}n
 +\f{x^2-1}x+\f{(x^2-4x-1)t^2}{(x+1)^2}.
 \end{align*}
 Now applying Theorem 2.1 yields the result.
 \vskip0.2cm

 \par{\bf Corollary 2.3} {\sl Let $\{L_n(t)\}$ be the Laguerre
 polynomials, $x\in\Bbb R$ and $x\not=0,-1$. Then}
 \begin{align*}&\sum_{n=0}^{p-1}\Big(\frac
 {(1-x)^3}{x(x+1)}n^2+\frac{4x^3-6x^2+2+4x(1-x^2)t}{x(x+1)^2}n
  \\&\q+\frac 1x+\frac{(x^2-1)(1-t)^2-4x(1-t)}{(x+1)^2}\Big)
  \frac{L_n(t)^2}{x^n}
 \\&=\frac 1{x^{p-1}}\Big(p^2L_p(t)^2+\frac{p^2}xL_{p-1}(t)^2
 -\Big(\frac 4{x+1}p^2+\frac{2(x-1)-2(x+1)t}{(x+1)^2}p\Big)
 L_p(t)L_{p-1}(t)\Big).
 \end{align*}

 \par{\it Proof.} Taking $c=1$, $r=1$, $P(n)=2n+1-t$, $a_0=1-t,\ a_1=2,\ a_2=a_3=0$, $u_n=L_n(t)$ in Theorem 2.1 gives
 $$s_n(x)=\f{4}{x+1}n+\f{2(x-1)-2(x+1)t}{(x+1)^2}$$
 and
 \begin{align*}&P_1(n,x)\\&=\frac{(n+1)^2}x-n^2x
 -(2n+1-t)^2+\Big(
 \frac{4xn}{x+1}
 +\frac{2x((x-1)-(x+1)t)}{(x+1)^2}\Big)(2n+1-t)
 \\&=\f{(1-x)^3}{x(x+1)}n^2+\f{4x^3-6x^2+2+4x(1-x^2)t}{x(x+1)^2}n
 +\f 1x+\f{(x^2-1)(1-t)^2-4x(1-t)}{(x+1)^2}.
 \\&\end{align*}
Now applying Theorem 2.1 yields the result.
\vskip0.2cm
\par{\bf Corollary 2.4} {\sl Let
$$d_{-1}(t)=0\quad\text{and}\quad d_n(t)=\sum_{k=0}^n\binom nk\binom tk2^k\quad(n=0,1,2,\ldots).$$
Then for every positive integer $p$ and $x\not=0,-1$,
\begin{align*}&\sum_{n=0}^{p-1}\Big(\frac{x^2-1}xn^2-\frac 2xn-\frac 1x+\frac{x-1}{x+1}(1+2t)^2\Big)\frac{d_n(t)^2}{(-x)^n}
\\&=\frac 1{(-x)^{p-1}}\Big(p^2d_p(t)^2-\frac{p^2}xd_{p-1}(t)^2
-\frac{2(1+2t)}{x+1}pd_p(t)
d_{p-1}(t)\Big).
\end{align*}}
\par{\it Proof.} From [26] and [19],
$$(n+1)d_{n+1}(t)=(1+2t)d_n(t)+nd_{n-1}(t)\q (n=0,1,2,\ldots).$$
Taking $c=-1$, $r=1$, $P(n)=1+2t$, $a_0=1+2t,\ a_1=a_2=a_3=0$, $u_n=d_n(t)$ in Theorem 2.1 gives
 $s_n(x)=\f{2(1+2t)}{x+1}$
 and
 \begin{align*}P_1(n,x)&=-\f 1x(n+1)^2+n^2x-(1+2t)^2+\f{2x(1+2t)^2}{x+1}
 \\&=\f{x^2-1}xn^2-\f 2xn-\f 1x+\f{x-1}{x+1}(1+2t)^2.
 \end{align*}
  Now applying Theorem 2.1 yields the result.
  \vskip0.2cm
  \par{\bf Corollary 2.5} {\sl Let $J_n(t)$ be the Bessel function given by
  $$J_n(t)=\sum_{k=0}^{\infty}\f{(-1)^k}{k!(n+k)!}\Ls t2^{n+2k}\ (n=0,1,2,\ldots).$$
  Then for every positive integer $p$ and $x\not=0,-1$,}
  \begin{align*}&\sum_{n=0}^{p-1}
  \Big(\f{4(x-1)}{(x+1)t^2}n^2-\f{8x}{(x+1)^2t^2}n+\f{1-x^2}x
  \Big)\f{J_n(t)^2}{x^n}
  \\&=\f 1{x^{p-1}}\Big(J_p(t)^2+\f 1xJ_{p-1}(t)^2-\f 4{(x+1)t}\Big(p-\f 1{x+1}\Big)J_p(t)J_{p-1}(t)\Big)
  \\&\q-xJ_0(t)^2-J_1(t)^2+\f {4x}{(x+1)^2t}J_0(t)J_1(t).
  \end{align*}
  \par{\it Proof.} Set $J_{-1}(t)=-J_1(t)$. Then $J_{n+1}(t)=\f{2n}tJ_n(t)-J_{n-1}(t)$ for $n\ge 0$.
  Taking $r=0,\ c=1,\ P(n)=\f{2n}t,\ a_0=a_2=a_3=0,\ a_1=\f 2t$ and $u_n=J_n(t)$ in Theorem 2.1 gives
\begin{align*}&s_n(x)=\f 4{(x+1)t}\Big(n-\f 1{x+1}\Big),
\\&P_r(n,x)=\f 1x-x-\f{4n^2}{t^2}+\f{4x}{(x+1)t}\Big(n-\f 1{x+1}\Big)\f{2n}t
\\&\qq\q\;=\f{4(x-1)}{(x+1)t^2}n^2-\f{8x}{(x+1)^2t^2}n+\f{1-x^2}x
\end{align*}
and so the result follows.

 \section*{3. Curious identities for Ap\'ery-like numbers}
 \par\q\ In this section, we use Theorem 2.1 to deduce the identities involving Ap\'ery-like numbers.
 \vskip0.2cm
 \par{\bf Theorem 3.1} {\sl Let $a,b,c,x\in\Bbb R$ with $cx(x+1)\not=0$ and $\{u_n\}$ be the sequence given by
 $$u_0=1,\ u_1=b\ \text{and}\  (n+1)^2u_{n+1}=(an(n+1)+b)u_n-cn^2u_{n-1}\ (n=1,2,3,\ldots).$$
 Let
  \begin{align*}
 P_n(a,b,c;x)
&= \frac{(x-1)\bigl(a^2x-c(x+1)^2\bigr)}{x(x+1)}n^4
+\frac{2a^2x(x^2-2x-1)+4c(x+1)^2}{x(x+1)^2}n^3\\
&\quad+\frac{(x^4-7x^3-5x^2-x)a^2+2x(x+1)(x^2-1)ab+6
(x+1)^3c}{x(x+1)^3}n^2\\
&\quad+\frac{-4a^2x^3+2x(x+1)(x^2-2x-1)ab+4(x+1)^3c      }{x(x+1)^3}n\\
&\quad+\frac{-4abx^3+x(x+1)(x^2-1)b^2+(x+1)^3c     }{x(x+1)^3}.
\end{align*}
Then for every positive integer $p$,
\begin{align*}
&\sum_{n=0}^{p-1}P_n(a,b,c;x)\f{u_n^2}{(cx)^n}
\\&=\f 1{(cx)^{p-1}}\Big(p^4u_p^2+\f cxp^4u_{p-1}^2-p^2
\Big(\frac{2a}{x+1}p^2
  + \frac{2a(x-1)}{(x+1)^2}p
 + \frac{2b(x+1)^2-4ax}{(x+1)^3}\Big)u_pu_{p-1}\Big).
 \end{align*}}
\par{\it Proof.} Putting $r=2$, $P(n)=an(n+1)+b$, $a_0=b,\ a_1=a_2=a,\ a_3=0$ in Theorem 2.1 gives
$$s_n(x)
= \frac{2a}{x+1}n^2
 + \frac{2a(x-1)}{(x+1)^2}n
 + \frac{2b(x+1)^2-4ax}{(x+1)^3}$$
and
$$P_2(n,x)=\f cx(n+1)^4-cxn^4-P(n)^2+xs_n(x)P(n)=P_n(a,b,c;x).$$
Now applying Theorem 2.1 yields the result.
\vskip0.2cm
\par {\bf Remark 3.1} If $\{u_n\}$ is given by (1.2) and $c=\f{a^2}4\not=0$, taking $x=1$ in Theorem 3.1 yields the identity:
$$\Big(c-\f{ab}2\Big)\sum_{n=0}^{p-1}(2n+1)\f{u_n^2}{c^n}
=\f 1{c^{p-1}}\Big(p^4u_p^2+cp^4u_{p-1}^2-p^2\Big(ap^2-\f a2+b\Big)u_pu_{p-1}\Big).$$
In the case $a=2$, $c=1$ and $b=\f{x+1}2$, the identity reduces to (1.6).
\vskip0.2cm
\par{\bf Corollary 3.1} {\sl For $p=1,2,3,\ldots,$
\begin{align*}&\sum_{n=0}^{p-1}(18n^3-21n^2-18n-4)\f{S_n^2}{64^n}
\\&=\f 1{16\cdot 64^{p-1}}\big(9p^4S_p^2+144p^4S_{p-1}^2-p^2
(72p^2+24p-8)S_pS_{p-1}\big)
\end{align*}
and
\begin{align*}&\sum_{n=0}^{p-1}(18n^3+75n^2+78n+25)
\frac{S_n^2}{16^n}
\\&=\frac 1{16^p}\big(9p^4S_p^2+576p^4S_{p-1}^2-p^2
(144p^2-48p-16)S_pS_{p-1}\big).
\end{align*}}
\par{\it Proof.} Taking $a=12,\ b=4,\ c=32$ and $x=2,\f 12$ in Theorem 3.1 gives the results.
\vskip0.2cm

\par{\bf Corollary 3.2} {\sl For $p=1,2,3,\ldots$ we have
\begin{align*}&\sum_{n=0}^{p-1}(7938n^3+8505n^2+4200n+813)
\f{f_n^2}{64^n}
\\&=\f{p^2}{64^{p-1}}\big(49p^2(f_p^2+f_{p-1}^2)
+(98p^2+126p+60)f_pf_{p-1}\big).\end{align*}}
\par{\it Proof.} Taking $a=7,\ b=2,\ c=-8$ and $x=-8$ in Theorem 3.1 gives the result.
\vskip0.2cm
\par{\bf Corollary 3.3} {\sl For $p=1,2,3,\ldots$ we have
\begin{align*}&\sum_{n=0}^{p-1}(1600n^3+820n^2+110n-19)
\frac{a_n^2}{81^n}
\\&=\frac{p^2}{81^{p-1}}\Big(\frac{25}2p^2(a_p^2+a_{p-1}^2)
-(25p^2+20p+3)a_pa_{p-1}\Big).\end{align*}}
\par{\it Proof.} Taking $a=10,\ b=3,\ c=9$ and $x=9$ in Theorem 3.1 gives the result.
\vskip0.2cm
\par{\bf Corollary 3.4} {\sl For $p=1,2,3,\ldots$ we have
\begin{align*}&\sum_{n=0}^{p-1}(578n^3+9673n^2+21624n+8973)
\frac{Q_n^2}{64^n}
\\&=\frac{p^2}{64^{p-1}}\Big(289p^2(Q_p^2+81Q_{p-1}^2)
+18(289p^2-17p-42)Q_pQ_{p-1}\Big)\end{align*}
and
\begin{align*}&\sum_{n=0}^{p-1}(578n^3-7939n^2+4012n+3556)
\frac{Q_n^2}{81^n}
\\&=\frac{p^2}{81^{p-1}}\Big(289p^2(Q_p^2+64Q_{p-1}^2)
+16(289p^2+17p-42)Q_pQ_{p-1}\Big).\end{align*}}
\par{\it Proof.} Taking $a=-17,\ b=-6,\ c=72$ and $x=\f 89,\f 98$ in Theorem 3.1 gives the results.
\vskip0.2cm

 \par{\bf Theorem 3.2} {\sl Let $a,b,c\in\Bbb R$ with $c\not=0$ and $\{u_n\}$ be the sequence given by
 $$u_0=1,\ u_1=b\ \text{and}\  (n+1)^2u_{n+1}=(an(n+1)+b)u_n-cn^2u_{n-1}\ (n=1,2,3,\ldots).$$
 Then for $p\in\Bbb Z^+$,
 \begin{align*} \sum_{n=0}^{p-1}Q(n)u_n^2
 &=(c+1)^3p^4u_p^2+c^2(c+1)^3p^4u_{p-1}^2
 \\&\q -p^2\big(2ac(c+1)^2p^2-2ac(c^2-1)p-4ac^2+2bc(c+1)^2\big)
 u_pu_{p-1},
 \end{align*}
    where
\begin{align*}
Q(n)=&(c+1)(c^2-1)((c+1)^2-a^2)n^4
\\&+\big(-2a^2(c+1)(c^2+2c-1)+4c^2(1+c)^3\big)n^3
\\&+\big(-a^2(c^3+5c^2+7c-1)-2(c+1)(c^2-1)ab+6c^2(c+1)^3\big)n^2
\\&+\big(-4a^2c-2(c+1)(c^2+2c-1)ab+4c^2(c+1)^3)n
\\&-(c+1)(c^2-1)b^2-4abc+c^2(c+1)^3.
\end{align*}}
\par{\it Proof.} Observe that
$Q(n)=(c+1)^3P_2\left(n,\frac 1c\right)$. Taking $x=\f 1c$ in Theorem 3.1 yields the result.
\vskip0.2cm
\par{\bf Corollary 3.5} For every positive integer $p$,
\begin{align*}&\sum_{n=0}^{p-1}(2646n^3+5103n^2+3668n+940)f_n^2
\\&=\frac {49}3p^4(f_p-8f_{p-1})^2+16p^2(21p-10)f_pf_{p-1},
\\&\sum_{n=0}^{p-1}\big(1600n^3+3980n^2+3270n+909\big)a_n^2
\\&=\frac{25}2p^4(a_p^2+81a_{p-1}^2)-9p^2(25p^2-20p+3)a_pa_{p-1},
\\&\sum_{n=0}^{p-1}
\left(10634085n^4+45622368n^3+70689360n^2
+47911968n+12084400\right)S_n^2
\\&=11979p^4(S_p^2+1024S_{p-1}^2)-256p^2(1089p^2-1023p+299)
S_pS_{p-1},
\\&\sum_{n=0}^{p-1}\left(3582488n^4+15116220n^3+23253129n^2+15705225n
+3954159\right)W_n^2
\\&=5488p^4(W_p^2+729W_{p-1}^2)+p^2(95256p^2
-88452p+25191)W_pW_{p-1},
\\&\sum_{n=0}^{p-1}\big(1095183237375n^4
+4415024793600n^3+6644229277440n^2
\\&\quad\quad+4436750566144n+1110021214352\big)G_n^2
\\&=16974593p^4\left(G_p^2+65536G_{p-1}^2\right)
-p^2\left(1082146816p^2-1073725440p+397416448\right)G_pG_{p-1}
\end{align*}
and
\begin{align*}&\sum_{n=0}^{p-1}\big(
1906929360n^4+7841889074n^3+11907294613n^2
\\&\qquad+7987243596n+2003013828\big)Q_n^2
\\&=389017p^4\left(Q_p^2+5184Q_{p-1}^2\right)
+144p^2\left(90593p^2-88111p+29526\right)Q_pQ_{p-1}.
\end{align*}
\par{\it Proof.} Putting $(a,b,c)=(7,2,-8),(10,3,9),(12,4,32),(-9,-3,27),(32,12,256),
(-17,$ $-6,72)$ in Theorem 3.2 gives the results.
\vskip0.2cm
\par{\bf Theorem 3.3} {\sl Let $a,b,c\in\Bbb R$ with $c\not=0$ and $\{u_n\}$ be the sequence given by
 $$u_0=1,\ u_1=b\ \text{and}\  (n+1)^2u_{n+1}=(an(n+1)+b)u_n-cn^2u_{n-1}\ (n=1,2,3,\ldots).$$
 Then for $p\in\Bbb Z^+$,
\begin{align*}
&(4c-a^2)\sum_{n=0}^{p-1}(2n+1)^3\frac{u_n^2}{c^n}
+(a^2-4ab+4c)\sum_{n=0}^{p-1}(2n+1)\frac{u_n^2}{c^n}
\\&=\frac {8p^2}{c^{p-1}}\Big(p^2u_p^2+cp^2u_{p-1}^2-\big(ap^2-\frac a2+b\big)u_pu_{p-1}\Big).
\end{align*}}
\par{\it Proof.}
Taking $x=1$ in Theorem 3.1 we see that
\begin{align*}P_n(a,b,c;1)&=\f 12\big((-2a^2+8c)n^3+(-3a^2+12c)n^2+(-a^2-2ab+8c)n-ab+2c\big)
\\&=\f 18\big((4c-a^2)(2n+1)^3+(a^2-4ab+4c)(2n+1)\big)
\end{align*}
and so the result follows.
\vskip0.2cm
\par{\bf Corollary 3.6} {\sl For every positive integer $p$,
\begin{align*}
&\sum_{n=0}^{p-1}\big(25(2n+1)^3+3(2n+1)\big)(-1)^n{A'_n}^2
\\&=(-1)^{p-1}\frac{4p^2}5\big(-2p^2{A'_p}^2+2p^2{A'_{p-1}}^2
+(22p^2-5)A'_pA'_{p-1}\big),
\\&\sum_{n=0}^{p-1}\big(27(2n+1)^3+13(2n+1)\big)\frac{f_n^2}{(-8)^n}
\\&=\frac{4p^2}{3\cdot (-8)^{p-1}}\big(-2p^2f_p^2+16p^2f_{p-1}^2
+(14p^2-3)f_pf_{p-1}\big),
\\&\sum_{n=0}^{p-1}\big(4(2n+1)^3-(2n+1)\big)
\frac{a_n^2}{9^n}
\\&=\frac{p^2}{2\cdot 9^{p-1}}\big(-p^2a_p^2-9p^2a_{p-1}^2
+(10p^2-2)a_pa_{p-1}\big),
\\&\sum_{n=0}^{p-1}\big((2n+1)^3-5(2n+1)\big)
\frac{S_n^2}{32^n}
\\&=\frac{p^2}{4\cdot 32^{p-1}}\big(-2p^2S_p^2-64p^2S_{p-1}^2
+(24p^2-4)S_pS_{p-1}\big),
\\&\sum_{n=0}^{p-1}\big((2n+1)^3+3(2n+1)\big)
\frac{W_n^2}{27^n}
\\&=\frac{4p^2}{27^{p}}\big(2p^2W_p^2+54p^2W_{p-1}^2
+(18p^2-3)W_pW_{p-1}\big)
\end{align*}
and
\begin{align*}&\sum_{n=0}^{p-1}\big((2n+1)^3-169(2n+1)\big)
\frac{Q_n^2}{72^n}
\\&=\frac{4p^2}{72^{p-1}}\big(-2p^2Q_p^2-144p^2Q_{p-1}^2
+(-34p^2+5)Q_pQ_{p-1}\big).
\end{align*}}
\par{\it Proof.} Taking $(a,b,c)=(11,3,-1),(7,2,-8),(10,3,9),(12,4,32),(32,12,256),
(-9,-3,27),$ $(-17,-6,72)$ in Theorem 3.3 gives the results.
\vskip0.2cm

\vskip0.2cm
\par{\bf Theorem 3.4} {\sl Let $a,b,c,x\in\Bbb R$ with $cx(x+1)\not=0$ and $\{u_n\}$ be the sequence given by
 $$u_0=1,\ u_1=b\ \text{and}\  (n+1)^3u_{n+1}=(2n+1)(an(n+1)+b)u_n-cn^3u_{n-1}\ (n\ge 1).$$
Let $p$ be a positive integer,
\begin{align*}
s_p(x)
&= \frac{4a}{x+1}p^3
 + \frac{6a(x-1)}{(x+1)^2}p^2
+ \frac{2a(x^2-10x+1)+4b(x+1)^2}{(x+1)^3}p
  \\&\quad+ \frac{2(x-1)\bigl(b(x+1)^2-6ax\bigr)}{(x+1)^4}
\end{align*}
and
$$ Q_n(a,b,c;x)=c_6n^6+c_5n^5+c_4n^4+c_3n^3+c_2n^2+c_1n+c_0,
$$ where
\begin{align*}
c_6&=\frac{(x-1)\bigl(4a^2x-c(x+1)^2\bigr)}{x(x+1)},
\quad c_5=\frac{12a^2x(x^2-2x-1)+6c(x+1)^2}{x(x+1)^2},
\\
c_4&=\frac 1{(x+1)^3}\big(a^2(13x^3-71x^2-49x-13)+8ab(x+1)^2(x-1)
+\frac{15c}x(x+1)^3\big),
\\c_3&=\frac 1{(x+1)^4}\big(6a^2(x^4-16x^3-12x^2-4x-1)+16ab(x+1)^2
(x^2-2x-1)
+\frac{20c}x(x+1)^4\big)
\\ c_2&=\frac 1{(x+1)^4}\big(a^2(x^4-58x^3+12x^2-2x-1)
+2ab(x+1)(5x^3-31x^2-17x-5)\\
&\q+4b^2(x-1)(x+1)^3+\frac{15c}x(x+1)^4\big),
\\ c_1&=\frac 1{(x+1)^4}\big(-12a^2x^2(x-1)+2ab(x^4-24x^3-4x^2-4x-1)\\
&\quad+4b^2(x+1)^2(x^2-2x-1)+\frac {6c}x(x+1)^4\big),
\\c_0&=\frac 1{(x+1)^4}\big(-12abx^2(x-1)+b^2(x+1)^2(x^2-4x-1) +\frac cx(x+1)^4
\big).
\end{align*}
Then
$$\sum_{n=0}^{p-1}Q_n(a,b,c;x)\f{u_n^2}{(cx)^n}
=\frac {p^3}{(cx)^{p-1}}
\Big(p^3u_p^2+\frac cxp^3u_{p-1}^2
-s_p(x)u_pu_{p-1}\Big).$$}
\par{\it Proof.} Taking $r=3$, $P(n)=(2n+1)(an(n+1)+b)=a_3n^3
+a_2n^2+a_1n+a_0$ in Theorem 2.1 and noting that $P_3(n,x)=Q_n(a,b,c;x)$ yields the result.
\vskip0.2cm
\par{\bf Corollary 3.7} {\sl For every positive integer $p$,
\begin{align*}&\sum_{n=0}^{p-1}\big(4500n^5
+225n^4-4570n^3-4215n^2-1588n
-224\big)\frac{D_n^2}{256^n}
\\&=\frac {p^3}{256^{p-1}}\Big(\frac{125}{12}p^3(D_p-4D_{p-1})^2
+(-75p^2+5p+14)D_pD_{p-1} \Big)
\end{align*}
and
\begin{align*}&\sum_{n=0}^{p-1}\big(4500n^5+22275n^4+39530n^3
+34155n^2+14732n+2556\big)\f{D_n^2}{16^n}
\\&=\f {p^3}{16^{p-1}}\Big(
\frac{125}{12}p^3(D_p-16D_{p-1})^2+(300p^2+20p-56)D_pD_{p-1}\Big).
\end{align*}}
\par{\it Proof.} Taking $a=10,\ b=4,\ c=64$, $x=4,\f 14$ and $u_n=D_n$ in Theorem 3.4 yields the results.

\vskip0.2cm
\par{\bf Theorem 3.5} {\sl Let $a,b,c\in\Bbb R$ with $c\not=0$ and $\{u_n\}$ be the sequence given by
 $$u_0=1,\ u_1=b\ \text{and}\  (n+1)^3u_{n+1}=(2n+1)(an(n+1)+b)u_n-cn^3u_{n-1}\ (n\ge 1).$$
For every positive integer $p$,
 \begin{align*}
&\sum_{n=0}^{p-1}
\big(3(c-a^2)(2n+1)^5+2(3a^2-8ab+5c)(2n+1)^3+(3(c-a^2)
+16(a-b)b)(2n+1)\big)\frac{u_n^2}{c^n}
\\&=\frac{16p^3}{c^{p-1}}\big(p^3(u_p^2+cu_{p-1}^2)
-(2ap^3-2(a-b)p)u_pu_{p-1}\big).
\end{align*}}
\par{\it Proof.} Putting $x=1$ in Theorem 3.4 gives
\begin{align*}Q_n(a,b,c;1)&=6(c-a^2)n^5
+15(c-a^2)n^4
+4(5c-3a^2-2ab)n^3\\
&\quad+3(5c-a^2-4ab)n^2
+2(3c-2ab-b^2)n
+(c-b^2)
\\&=\frac 1{16}\big(3(c-a^2)(2n+1)^5+2(3a^2-8ab+5c)(2n+1)^3+(3(c-a^2)
\\&\quad+16(a-b)b)(2n+1)\big).
\end{align*}
Thus the result follows from Theorem 3.4.
\vskip0.2cm
\par{\bf Remark 3.2} Taking $a=c=1$ and $b=x$ in Theorem 3.5 yields (1.7).
\vskip0.2cm
\par{\bf Corollary 3.8} {\sl For every positive integer $p$,
\begin{align*}&\sum_{n=0}^{p-1}\big(9(2n+1)^5-4(2n+1)^3-
(2n+1)\big)A_n^2
\\&=-\frac{p^4}6\big(p^2(A_p^2+A_{p-1}^2)-(34p^2-24)
A_pA_{p-1}\big),
\\&\sum_{n=0}^{p-1}\big(9(2n+1)^5-50(2n+1)^3-23
(2n+1)\big)\frac{D_n^2}{64^n}
\\&=-\frac{4p^4}{3\cdot 64^{p-1}}\big(p^2(D_p^2+64D_{p-1}^2)-(20p^2-12)
D_pD_{p-1}\big),
\\&\sum_{n=0}^{p-1}\big((2n+1)^5+8(2n+1)^3
+3(2n+1)\big)\frac{b_n^2}{81^n}
\\&=\frac{p^4}{6\cdot 81^{p-1}}\big(p^2(b _p^2+81b_{p-1}^2)+(14p^2-8)
b_pb_{p-1}\big),
\\&\sum_{n=0}^{p-1}\big(3(2n+1)^5-2(2n+1)^3
-(2n+1)\big)\frac{T_n^2}{16^n}
\\&=-\frac{p^4}{8\cdot 16^{p-1}}\big(p^2(T_p^2+16T_{p-1}^2)-(24p^2-16)
T_pT_{p-1}\big).
\end{align*}}
\par{\it Proof.} Taking $(a,b,c)=(17,5,1),(10,4,64),
(-7,-3,81),(12,4,16)$ in Theorem 3.5 yields the results.
\vskip0.2cm

\par{\bf Theorem 3.6}  {\sl Let $a,b,c\in\Bbb R$ with $c\not=0$ and $\{u_n\}$ be the sequence given by
 $$u_0=1,\ u_1=b\ \text{and}\  (n+1)^3u_{n+1}=(2n+1)(an(n+1)+b)u_n-cn^3u_{n-1}\ (n\ge 1).$$
Let
\begin{align*}
t_p(a,b,c)&=4ac(1+c)^3p^3
+6ac(1-c)(1+c)^2p^2
\\&\quad+2c(1+c)\bigl(a(1-10c+c^2)+2b(1+c)^2\bigr)p
+2c(1-c)\bigl(b(1+c)^2-6ac\bigr)
\end{align*}
and
$$T_n(a,b,c)
=C_6n^6+C_5n^5+C_4n^4+C_3n^3+C_2n^2+C_1n+C_0,$$
where
\begin{align*}
C_6&=(1+c)^3(1-c)\left(4a^2-(1+c)^2\right),
\quad \\ C_5&=(1+c)^2\big(6\left(2a^2(1-2c-c^2)+c^2(1+c)^2\right)\big),
\\C_4&=(1+c)\big(a^2(13-71c-49c^2-13c^3)
+8ab(1+c)^2(1-c)
+15c^2(1+c)^3\big),
\\ C_3&=6a^2(1-16c-12c^2-4c^3-c^4)+16ab(1+c)^2(1-2c-c^2)
+20c^2(1+c)^4,
\\ C_2&=a^2(1-58c+12c^2-2c^3-c^4)
+2ab(1+c)(5-31c-17c^2-5c^3)
\\&\quad+4b^2(1-c)(1+c)^3
+15c^2(1+c)^4,
\\C_1&=-12a^2c(1-c)
+2ab(1-24c-4c^2-4c^3-c^4)
+4b^2(1+c)^2(1-2c-c^2)
+6c^2(1+c)^4,
\\ C_0&=-12abc(1-c)
+b^2(1+c)^2(1-4c-c^2)
+c^2(1+c)^4.
\end{align*}
Then for $p\in\Bbb Z^+$,
$$\sum_{n=0}^{p-1}T_n(a,b,c)u_n^2
=p^3\big((1+c)^4p^3(u_p^2+c^2u_{p-1}^2)
-t_p(a,b,c)u_pu_{p-1}\big).$$}
\par{\it Proof.} One can check that
$$t_p(a,b,c)=(1+c)^4s_p\Big(a,b,c;\f 1c\Big)
\q\t{and}\q T_n(a,b,c)=(1+c)^4Q_n\Big(a,b,c;\f 1c\Big).$$
Thus the result follows from Theorem 3.4.
\vskip0.2cm
\par{\bf Corollary 3.9} {\sl For every positive integer $p$,
\begin{align*}
&\sum_{n=0}^{p-1}\big(
13\,235\,551\,875n^6
+83\,457\,270\,000n^5
+213\,544\,152\,900n^4
\\&\qquad
+288\,035\,725\,560n^3
+217\,344\,276\,900n^2
+87\,226\,480\,464n
+14\,564\,793\,552
\big)D_n^2
\\&=p^3\big(
3\,570\,125p^3\big(D_p^2+4096D_{p-1}^2\big)
\\ &\qquad-\big(140\,608\,000p^3
-204\,422\,400p^2 +113\,767\,680p
-21\,063\,168 \big)D_pD_{p-1}\big),
\\&\sum_{n=0}^{p-1}
\big(-21\,150\,465n^6 -15\,037\,248n^5
+128\,613\,840n^4 +290\,273\,888n^3
\\&\qquad
+264\,369\,072n^2 +115\,395\,232n
+20\,044\,560\big)T_n^2
\\&=p^3\big(83\,521p^3\big(T_p^2+256T_{p-1}^2\big)
-128(29\,478p^3-39\,015p^2+14\,773p-15)T_pT_{p-1}\big),
\\&\sum_{n=0}^{p-1} \big( 4\,499\,162\,880 n^6 + 27\,394\,463\,568 n^5 + 68\,954\,247\,310 n^4
\\&\quad+ 92\,253\,917\,556 n^3 + 69\,316\,663\,129 n^2 + 27\,754\,699\,191 n + 4\,628\,470\,743 \big) b_n^2
\\&= \frac{p^3}{4} \big( 2\,825\,761 p^3 (b_p^2+6561 b_{p-1}^2)
\\&\quad + \big( 78\,156\,414 p^3 - 114\,375\,240 p^2 + 66\,924\,792 p - 13\,583\,700 \big) b_p b_{p-1} \big)
\end{align*}
and
\begin{align*}&\sum_{n=0}^{p-1}
\big(281462092005375n^6
+1701992057932800n^5
+4269488541077760n^4
\\&\qquad+5702341771922944n^3
+4280674043102976n^2
+1713159842955776n
\\&\qquad+285617606524992\big)V_n^2
\\&=p^3\big(4362470401p^3\big(V_p^2+65536V_{p-1}^2\big)
\\&\qquad-\big(278111731712p^3
-413921157120p^2+271644114944p-65778216960
\big)V_pV_{p-1}\big).
\end{align*}}
\par{\it Proof.} Taking $(a,b,c)=(10,4,64),(12,4,16),
(-7,-3,81),(16,8,256)$ in Theorem 3.6 yields the results.
\vskip0.2cm
\section*{4. Congruences involving Ap\'ery-like numbers}
\par\q\  Using the identities in Section 3 and known congruences for Ap\'ery-like numbers, in this section we deduce new supercongruences involving Ap\'ery-like numbers.
\vskip0.2cm
\par {\bf Lemma 4.1 ([8])} {\sl
Let $p > 3$ be a prime. Then}
$$A_{p-1} \equiv 1 + \frac{2}{3} p^3 B_{p-3} \pmod{p^4}
\q\t{and}\q
A_p \equiv 5 - \frac{14}{3} p^3 B_{p-3} \pmod{p^4}.$$
\par{\bf Theorem 4.1} {\sl Let $p$ be a prime greater than $3$. Then
$$\sum_{n=0}^{p-1}\big(9(2n+1)^5-4(2n+1)^3-
(2n+1)\big)A_n^2\e -20p^4+24p^6+\f{16}3p^7B_{p-3}\mod {p^8}.$$}
\par{\it Proof.}
From Corollary 3.8 and Lemma 4.1,
\begin{align*}&\sum_{n=0}^{p-1}\big(9(2n+1)^5-4(2n+1)^3-
(2n+1)\big)A_n^2
\\&=-\frac{p^4}6\big(p^2(A_p^2+A_{p-1}^2)-(34p^2-24)
A_pA_{p-1}\big)
\\&\e -\f {p^4}6\Big(p^2(5^2+1^2)-(34p^2-24)\Big(1 + \frac{2}{3} p^3 B_{p-3}\Big)\Big(5 - \frac{14}{3} p^3 B_{p-3}\Big)\Big)
\mod {p^8},
\end{align*}
which yields the result.
\vskip0.2cm
\par{\bf Lemma 4.2 ([16,28])}
{\sl Let $p > 3$ be a prime. Then}
$$D_{p-1} \equiv 64^{p-1} - \frac{p^3}{6} B_{p-3} \pmod{p^4}
\q\t{and}\q D_p \equiv 4 + \frac{16}{3} p^3 B_{p-3} \pmod{p^4}.$$
\par{\bf Theorem 4.2} {\sl Let $p$ be a prime greater than $3$. Then}
\begin{align*}&\sum_{n=0}^{p-1}\big(9(2n+1)^5-50(2n+1)^3-23
(2n+1)\big)\frac{D_n^2}{64^n}
\\&\e -64p^4-384p^6(2^{p-1}-1)-\f {224}3p^7B_{p-3}\mod {p^8},
\\&\sum_{n=0}^{p-1}\big(4500n^5+22275n^4+39530n^3
+34155n^2+14732n+2556\big)\f{D_n^2}{16^n}
\\&\e 16\cdot 4^{p-1}(75p^2+5p-14)p^3+ 1500p^6-\f{784}3p^6B_{p-3}\mod {p^7}
\\&\sum_{n=0}^{p-1}\big(4500n^5
+225n^4-4570n^3-4215n^2-1588n
-224\big)\frac{D_n^2}{256^n}
\\&\e 4^{2-p}(-75p^2+5p+14)p^3+\f{196}3p^6B_{p-3}\mod {p^7}.
\end{align*}
\par{\it Proof.} From Corollary 3.8 and Lemma 4.2,
\begin{align*}&\sum_{n=0}^{p-1}\big(9(2n+1)^5-50(2n+1)^3-23
(2n+1)\big)\frac{D_n^2}{64^n}
\\&=-\frac{4p^4}{3\cdot 64^{p-1}}\big(p^2(D_p^2+64D_{p-1}^2)-(20p^2-12)
D_pD_{p-1}\big)
\\&\e -\frac{4p^4}{3\cdot 64^{p-1}}\Big(p^2(16+ 64^{2p-1})-(20p^2-12)\Big(4 + \frac{16}{3} p^3 B_{p-3}
\Big)\Big(64^{p-1} - \frac{p^3}{6} B_{p-3}\Big)\Big)\mod {p^8},
\end{align*}
which yields the first congruence. The remaining congruences follow from Corollary 3.7 and Lemma 4.2.
\vskip0.2cm
\par{\bf Lemma 4.3 ([6,16])} {\sl
Let $p > 3$ be a prime. Then}
$$T_{p-1} \equiv 16^{p-1} + \frac{p^3}{4} B_{p-3} \pmod{p^4}
\q\t{and}\q
T_p \equiv 4 - p^3 B_{p-3} \pmod{p^4}.$$
\par{\bf Theorem 4.3} {\sl Let $p>3$ be a prime. Then
$$\sum_{n=0}^{p-1}\big(3(2n+1)^5-2(2n+1)^3
-(2n+1)\big)\frac{T_n^2}{16^n}
\e -8p^4+8p^6\mod {p^8}.$$}
\par{\it Proof.} From Corollary 3.8, Lemma 4.3 and the fact that $16^{p-1}\e 1+4(2^{p-1}-1)\mod {p^2}$ we deduce that
\begin{align*}
&\sum_{n=0}^{p-1}\big(3(2n+1)^5-2(2n+1)^3
-(2n+1)\big)\frac{T_n^2}{16^n}
\\&=-\frac{p^4}{8\cdot 16^{p-1}}\big(p^2(T_p^2+16T_{p-1}^2)-(24p^2-16)
T_pT_{p-1}\big)
\\&\e -\frac{p^6}{8\cdot 16^{p-1}}(4^2+16\cdot 16^{2(p-1)})+\f{p^4}{16^{p-1}}(3p^2-2)\cdot 4\Big(1-\f{p^3}4B_{p-3}\Big)
\Big(16^{p-1}+\f{p^3}4B_{p-3}\Big)
\\&\e -2p^6(1-4(2^{p-1}-1))-2p^6(1+4(2^{p-1}-1))+4p^4(3p^2-2)
\\&=8p^6-8p^4\mod {p^8}.
\end{align*}
This proves the theorem.
\vskip0.2cm
\par{\bf Lemma 4.4 ([9])} {\sl
Let $p > 3$ be a prime. Then}
\begin{align*}
&b_{p-1} \equiv 81^{p-1} - \frac{2}{27}p^3 B_{p-3} \pmod{p^4}, \\&b_p \equiv -3 - 6p^3 B_{p-3} \pmod{p^4}.
\end{align*}
\par{\bf Theorem 4.4} {\sl Let $p>3$ be a prime. Then
\begin{align*}
&\sum_{n=0}^{p-1}\big((2n+1)^5+8(2n+1)^3
+3(2n+1)\big)\frac{b_n^2}{81^n}
\\&\e 4p^4+8p^6+48p^6(3^{p-1}-1)+\f{208}{27}p^7B_{p-3}\mod {p^8}.
\end{align*}}
\par{\it Proof.} From Corollary 3.8 and Lemma 4.4,
 \begin{align*}
&\sum_{n=0}^{p-1}\big((2n+1)^5+8(2n+1)^3
+3(2n+1)\big)\frac{b_n^2}{81^n}
\\&=\frac{p^4}{6\cdot 81^{p-1}}\big(p^2(b _p^2+81b_{p-1}^2)+(14p^2-8)
b_pb_{p-1}\big)
\\&\e \f{p^6}{6\cdot 81^{p-1}}\big(9+81\cdot 81^{2(p-1)}\big)+\f{p^4}{6\cdot 81^{p-1}}(14p^2-8)(-3-6p^3B_{p-3})\Big(81^{p-1}-\f 2{27}p^3B_{p-3}\Big)
\\&\e \f 32p^6(1-4(3^{p-1}-1))+\f {27}2p^6(1+4(3^{p-1}-1))
-p^4(7p^2-4)\Big(1+\f{52}{27}p^3B_{p-3}\Big)\mod {p^8},
\end{align*}
which yields the result.
\vskip0.2cm
\par{\bf Lemma 4.5 ([8])} {\sl
Let $p > 3$ be a prime. Then
$$A'_{p-1}\e 1+\f 53p^3B_{p-3}\mod {p^4}
\q\t{and}\q A'_p\e 3-\f 53p^3B_{p-3}\mod {p^4}.$$}

\par{\bf Theorem 4.5} {\sl Let $p>3$ be a prime. Then}
$$\sum_{n=0}^{p-1}\big(25(2n+1)^3+3(2n+1)\big)(-1)^n{A'_n}^2
\e -12p^2+40p^4-\f{40}3p^5B_{p-3}\mod {p^6}.$$
\par{\it Proof.} From Corollary 3.6 and Lemma 4.5,
\begin{align*}
&\sum_{n=0}^{p-1}\big(25(2n+1)^3+3(2n+1)\big)(-1)^n{A'_n}^2
\\&=\frac{4p^2}5\big(-2p^2{A'_p}^2+2p^2{A'_{p-1}}^2
+(22p^2-5)A'_pA'_{p-1}\big)
\\&\e -\f 85p^4(3^2-1^2)+\f 45p^2(22p^2-5)\Big(3-\f 53p^3B_{p-3}\Big)\Big(1+\f 53p^3B_{p-3}\Big)
\\&\e -\f{64}5p^4+\f 45p^2(22p^2-5)\Big(3+\f{10}3p^3B_{p-3}\Big)
\\&\e -12p^2+40p^4-\f{40}3p^5B_{p-3}\mod {p^6}.
\end{align*}
The proof is complete.
\vskip0.2cm
\par{\bf Lemma 4.6 ([6,17])} {\sl
Let $p > 3$ be a prime. Then
$$f_{p-1}\e 8^{p-1}+\f 58p^3B_{p-3}\mod {p^4}
\q\t{and}\q f_p\e 2+\f 12p^3B_{p-3}\mod {p^4}.$$}

\par{\bf Theorem 4.6} {\sl Let $p>3$ be a prime. Then}
\begin{align*}
&\sum_{n=0}^{p-1}\big(27(2n+1)^3+13(2n+1)\big)\f{f_n^2}{(-8)^n}
\\&\e -8p^2+48p^4+96p^4(2^{p-1}-1)-7p^5B_{p-3}\mod {p^6},
\\&\sum_{n=0}^{p-1}(2646n^3+5103n^2+3668n+940)f_n^2
\\&\e 32p^2(21p-10)8^{p-1}+196p^4(3+8(8^{p-1}-1))-280p^5B_{p-3} \mod {p^6},
\\&\sum_{n=0}^{p-1}(7938n^3+8505n^2+4200n+813)
\f{f_n^2}{64^n}
\\&\e 245p^4-392p^4(8^{p-1}-1)+4p^2(49p^2+63p+30)8^{1-p}
+105p^5B_{p-3}\mod {p^6}.
\end{align*}
\par{\it Proof.} From Corollary 3.6 and Lemma 4.6,
\begin{align*}
&\sum_{n=0}^{p-1}\big(27(2n+1)^3+13(2n+1)\big)\frac{f_n^2}{(-8)^n}
\\&=\frac{4p^2}{3\cdot (-8)^{p-1}}\big(-2p^2f_p^2+16p^2f_{p-1}^2
+(14p^2-3)f_pf_{p-1}\big)
\\&\e -\f{8p^4}{3\cdot 8^{p-1}}(2^2-8\cdot 8^{2(p-1)})+\f{4p^2}{3\cdot 8^{p-1}}(14p^2-3)
\Big(2+\f 12p^3B_{p-3}\Big)\Big(8^{p-1}+\f 58p^3B_{p-3}\Big)
\\&\e -\f{32}3p^4(1-3(2^{p-1}-1))+\f{64}3p^4(1+3(2^{p-1}-1))
+\f{4p^2}{3}(14p^2-3)\Big(2+\f 74p^3B_{p-3}\Big)
\\&\e -8p^2+48p^4+96p^4(2^{p-1}-1)-7p^5B_{p-3}\mod {p^6}.
\end{align*}
By Corollary 3.5 and Lemma 4.6,
\begin{align*}&\sum_{n=0}^{p-1}(2646n^3+5103n^2+3668n+940)f_n^2
\\&=\frac {49}3p^4(f_p-8f_{p-1})^2+16p^2(21p-10)f_pf_{p-1}
\\&\e \frac {49}3p^4(2-8\cdot 8^{p-1})^2+16p^2(21p-10)\Big(2\cdot 8^{p-1}+\f 74p^3B_{p-3}\Big)
\\&\e 32p^2(21p-10)8^{p-1}+196p^4(3+8(8^{p-1}-1))-280p^5B_{p-3} \mod {p^6}.
\end{align*}
Similarly, from Corollary 3.2 and Lemma 4.6 we deduce the remaining congruence.

\vskip0.2cm
\par{\bf Lemma 4.7 ([16])} {\sl
Let $p > 3$ be a prime. Then
$$S_{p-1}\e (-1)^{\f{p-1}2}32^{p-1}+p^2E_{p-3}\mod {p^3}
\q\t{and}\q S_p\e 4+8(-1)^{\f{p-1}2}p^2E_{p-3}
\mod {p^3}.$$}

\par{\bf Theorem 4.7} {\sl Let $p>3$ be a prime. Then
\begin{align*}
&\sum_{n=0}^{p-1}\big((2n+1)^3-5(2n+1)\big)
\frac{S_n^2}{32^n}
\\&\e -4(-1)^{\f{p-1}2}p^2
+\big(24((-1)^{\f{p-1}2}-1)-12E_{p-3}\big)p^4\mod {p^5},
\\& \sum_{n=0}^{p-1}(18n^3-21n^2-18n-4)\f{S_n^2}{64^n}
\\&\e 18p^4-(-1)^{\f{p-1}2}p^2(18p^2+6p-2)2^{1-p}+6p^4E_{p-3}
\mod {p^5},
\\&\sum_{n=0}^{p-1}(18n^3+75n^2+78n+25)
\frac{S_n^2}{16^n}
\\&\e  45p^4-4(-1)^{\f{p-1}2}p^2(9p^2-3p-1)2^{p-1}+12p^4E_{p-3}
\mod {p^5}
\end{align*}
and
\begin{align*}&\sum_{n=0}^{p-1}
\left(10634085n^4+45622368n^3+70689360n^2
+47911968n+12084400\right)S_n^2
\\&\e -1024(-1)^{\f{p-1}2}p^2(1089p^2-1023p+299)32^{p-1}
+(12458160-918528E_{p-3})p^4\mod {p^5}.
\end{align*}}
\par{\it Proof.} By Lemma 4.7,
\begin{align*}S_pS_{p-1}&\e 4\big(1+2(-1)^{\f{p-1}2}p^2E_{p-3}\big)\big((-1)^{\f{p-1}2}32^{p-1}+p^2E_{p-3}\big)
\\&\e 4(-1)^{\f{p-1}2}32^{p-1}+12p^2E_{p-3}\mod {p^3}.
\end{align*}

From Corollary 3.6 and Lemma 4.7,
\begin{align*}
&\sum_{n=0}^{p-1}\big((2n+1)^3-5(2n+1)\big)
\frac{S_n^2}{32^n}
\\&=\frac{p^2}{4\cdot 32^{p-1}}\big(-2p^2S_p^2-64p^2S_{p-1}^2
+(24p^2-4)S_pS_{p-1}\big)
\\&\e\f{p^4}{4\cdot 32^{p-1}}(-2\cdot 4^2-64\cdot
 32^{2(p-1)})
+\f{p^2}{4\cdot 32^{p-1}}(24p^2-4)\cdot \big(4(-1)^{\f{p-1}2}32^{p-1}+12p^2E_{p-3}\big)
\\&\e -24p^4+p^2(24p^2-4)\big((-1)^{\f{p-1}2}+3p^2E_{p-3}\big)
\\&=-4(-1)^{\f{p-1}2}p^2
+\big(24((-1)^{\f{p-1}2}-1)-12E_{p-3}\big)p^4\mod {p^5},
\end{align*}
By Corollary 3.1 and Lemma 4.7,
\begin{align*}&\sum_{n=0}^{p-1}(18n^3-21n^2-18n-4)\f{S_n^2}{64^n}
\\&=\f 1{16\cdot 64^{p-1}}\big(9p^4S_p^2+144p^4S_{p-1}^2-p^2
(72p^2+24p-8)S_pS_{p-1}\big)
\\&\e \f {p^4}{16\cdot 64^{p-1}}(9\cdot 4^2+144\cdot 32^{2(p-1)})
\\&\q-\f {p^2}{16\cdot 64^{p-1}}(72p^2+24p-8)\big(4(-1)^{\f{p-1}2}
32^{p-1}+12p^2E_{p-3}\big)
\\&\e 18p^4-(-1)^{\f{p-1}2}p^2(18p^2+6p-2)2^{1-p}+6p^4E_{p-3}
\mod {p^5}
\end{align*}
and
\begin{align*}&\sum_{n=0}^{p-1}(18n^3+75n^2+78n+25)
\frac{S_n^2}{16^n}
\\&=\frac 1{16^p}\big(9p^4S_p^2+576p^4S_{p-1}^2-p^2
(144p^2-48p-16)S_pS_{p-1}\big)
\\&\e \f{p^4}{16^p}(9\cdot 4^2+576\cdot 32^{2(p-1)})
-\f{p^2}{16^{p-1}}(9p^2-3p-1)\big(4(-1)^{\f{p-1}2}
32^{p-1}+12p^2E_{p-3}\big)
\\&\e 45p^4-4(-1)^{\f{p-1}2}p^2(9p^2-3p-1)2^{p-1}+12p^4E_{p-3}
\mod {p^5}.
\end{align*}
From Corollary 3.5 and Lemma 4.7,
\begin{align*}
&\sum_{n=0}^{p-1}
\left(10634085n^4+45622368n^3+70689360n^2
+47911968n+12084400\right)S_n^2
\\&=11979p^4(S_p^2+1024S_{p-1}^2)-256p^2(1089p^2-1023p+299)
S_pS_{p-1}
\\&\e 11979p^4(4^2+1024\cdot 32^{2(p-1)})-256p^2
(1089p^2-1023p+299)\\&\q\times 4((-1)^{\f{p-1}2}32^{p-1}+3p^2E_{p-3})
\\&\e 12458160p^4-1024(-1)^{\f{p-1}2}p^2(1089p^2-1023p+299)32^{p-1}
-918528p^4E_{p-3}\mod{p^5}.
\end{align*}
The proof is now complete.
\vskip0.2cm
\par From now on let $\{U_n\}$ be given by
$$U_{2n-1}=0,\q U_0=1\q{and}\q U_{2n}=-2\sum_{k=0}^{n-1}
\b {2n}{2k}U_{2k}\q(n\ge 1).$$
\par{\bf Lemma 4.8 ([7,12])} {\sl
Let $p > 3$ be a prime. Then
$$a_p\e 3+3(-1)^{[\f p3]}p^2U_{p-3}\mod {p^3}
\q\t{and}\q a_{p-1}\e (-1)^{[\f p3]}9^{p-1}+p^2U_{p-3}
\mod {p^3}.$$}
\par{\bf Theorem 4.8} {\sl
Let $p > 3$ be a prime. Then
\begin{align*}&\sum_{n=0}^{p-1}\big(4(2n+1)^3-(2n+1)\big)
\frac{a_n^2}{9^n}
\\&\e -3(-1)^{[\f p3]}p^2+\big(15(-1)^{[\f p3]}-9-6U_{p-3}\big)p^4\mod {p^5},
\\&\sum_{n=0}^{p-1}\big(1600n^3+3980n^2+3270n+909\big)a_n^2
\\&\e 1125p^4-27(-1)^{[\f p3]}p^2(25p^2-20p+3)9^{p-1}-162p^4U_{p-3}\mod {p^5},
\\&\sum_{n=0}^{p-1}(1600n^3+820n^2+110n-19)
\frac{a_n^2}{81^n}
\\&\e 125p^4-3\s3 p^2(25p^2+20p+3)9^{1-p}-18p^4U_{p-3}\mod {p^5}.
\end{align*}}
\par{\it Proof.} Using Lemma 4.8,
\begin{align*}a_pa_{p-1}&\e \big(3+3(-1)^{[\f p3]}p^2U_{p-3}\big)
\big( (-1)^{[\f p3]}9^{p-1}+p^2U_{p-3}\big)
\\&\e 3\big((-1)^{[\f p3]}9^{p-1}+2p^2U_{p-3}\big)\mod {p^3}.\end{align*}
From Corollary 3.6 and Lemma 4.8,
\begin{align*}&\sum_{n=0}^{p-1}\big(4(2n+1)^3-(2n+1)\big)
\frac{a_n^2}{9^n}
\\&=\frac{p^2}{2\cdot 9^{p-1}}\big(-p^2a_p^2-9p^2a_{p-1}^2
+(10p^2-2)a_pa_{p-1}\big)
\\&\e -\f{p^4}{2\cdot 9^{p-1}}(3^2+9\cdot 9^{2(p-1)})+\f{p^2}{9^{p-1}}(5p^2-1)\cdot 3((-1)^{[\f p3]}9^{p-1}+2p^2U_{p-3})
\\&\e -3(-1)^{[\f p3]}p^2+\big(15(-1)^{[\f p3]}-9-6U_{p-3}\big)p^4\mod {p^5}.
\end{align*}
By Corollary 3.5 and Lemma 4.8,
\begin{align*}
&\sum_{n=0}^{p-1}\big(1600n^3+3980n^2+3270n+909\big)a_n^2
\\&=\frac{25}2p^4(a_p^2+81a_{p-1}^2)
-9p^2(25p^2-20p+3)a_pa_{p-1}
\\&\e \f{25}2p^4\big(3^2+81\cdot 9^{2(p-1)}\big)
-9p^2(25p^2-20p+3)\cdot 3((-1)^{[\f p3]}9^{p-1}+2p^2U_{p-3})
\\&\e 1125p^4-27(-1)^{[\f p3]}p^2(25p^2-20p+3)9^{p-1}-162p^4U_{p-3}\mod {p^5}.
\end{align*}
Also, from Corollary 3.3 and Lemma 4.8,
\begin{align*}&\sum_{n=0}^{p-1}(1600n^3+820n^2+110n-19)
\frac{a_n^2}{81^n}
\\&=\frac{p^2}{81^{p-1}}\Big(\frac{25}2p^2(a_p^2+a_{p-1}^2)
-(25p^2+20p+3)a_pa_{p-1}\Big)
\\&\e\f{25p^4}{2\cdot 81^{p-1}}(3^2+9^{2(p-1)})-
\f{p^2}{81^{p-1}}(25p^2+20p+3)\cdot 3((-1)^{[\f p3]}9^{p-1}+2p^2U_{p-3})
\\&\e 125p^4-3\s3 p^2(25p^2+20p+3)9^{1-p}-18p^4U_{p-3}\mod {p^5}.
\end{align*}
This completes the proof.
\vskip0.2cm
\par{\bf Lemma 4.9 ([13])} {\sl
Let $p > 3$ be a prime. Then
$$W_{p-1}\e \s3 27^{p-1}+p^2U_{p-3}\mod {p^3}
\q\t{and}\q W_{p}\e -3-9\s3 p^2U_{p-3}\mod {p^3}.$$}
 \par{\bf Theorem 4.9 } {\sl
Let $p > 3$ be a prime. Then}
\begin{align*}&\sum_{n=0}^{p-1}\big((2n+1)^3+3(2n+1)\big)
\frac{W_n^2}{27^n}
\\&\e  \f 43\s3 p^2+\f{p^4}3\big(32-24\s3+16U_{p-3}\big)\mod {p^5},
\\&\sum_{n=0}^{p-1}\left(3582488n^4+15116220n^3+23253129n^2+15705225n
+3954159\right)W_n^2
\\&\e -3(95256p^4
-88452p^3+25191p^2)\s3 27^{p-1}+(4050144-302292U_{p-3})p^4\mod {p^5}.
\end{align*}
\par{\it Proof.} By Lemma 4.9,
\begin{align*}W_pW_{p-1}&\e -3\big(1+3\s3 p^2U_{p-3}\big)\big(\s3 27^{p-1}+p^2U_{p-3}\big)
\\&\e -3\big(\s3 27^{p-1}+4p^2U_{p-3}\big)\mod {p^3}.
\end{align*}
From Corollary 3.6 and Lemma 4.9,
\begin{align*}&\sum_{n=0}^{p-1}\big((2n+1)^3+3(2n+1)\big)
\frac{W_n^2}{27^n}
\\&=\frac{4p^2}{27^{p}}\big(2p^2W_p^2+54p^2W_{p-1}^2
+(18p^2-3)W_pW_{p-1}\big)
\\&\e \f{8p^4}{27}((-3)^2+27\cdot 27^{2(p-1)})+\f{4p^2}{27^p}(18p^2-3)\cdot (-3)
\big(\s3 27^{p-1}+4p^2U_{p-3}\big)
\\&\e \f{32}3p^4-\f 43\s3 p^2(6p^2-1)+\f{16}3p^4U_{p-3}\mod {p^5}.
\end{align*}
Also, from Corollary 3.5 and Lemma 4.9 we have
\begin{align*}&\sum_{n=0}^{p-1}\left(3582488n^4+15116220n^3+23253129n^2+15705225n
+3954159\right)W_n^2
\\&=5488p^4(W_p^2+729W_{p-1}^2)+p^2(95256p^2
-88452p+25191)W_pW_{p-1}
\\&\e 5488p^4((-3)^2+729\cdot27^{2(p-1)})+p^2(95256p^2
-88452p+25191)(-3)\big(\s3 27^{p-1}+4p^2U_{p-3}\big)
\\&\e 4050144p^4-3p^2(95256p^2
-88452p+25191)\s3 27^{p-1}-302292p^4U_{p-3}\mod {p^5}.
\end{align*}
\par{\bf Lemma 4.10 ([10,21])} {\sl Let $p$ be a prime greater than $3$. Then
$$G_p\e 12+64(-1)^{\f{p-1}2}p^2E_{p-3}\mod {p^3}
\ \t{and}\  G_{p-1}\e (-1)^{\f{p-1}2}256^{p-1}+3p^2E_{p-3}
\mod {p^3}.$$ }

\par{\bf Theorem 4.10} {\sl Let $p$ be a prime greater than $3$. Then}
\begin{align*}&\sum_{n=0}^{p-1}\big(1095183237375n^4
+4415024793600n^3+6644229277440n^2
\\&\quad\quad+4436750566144n+1110021214352\big)G_n^2
\\&\e -12\left(1082146816p^4-1073725440p^3
+397416448p^2\right)(-1)^{\f{p-1}2}256^{p-1}
\\&\q+(1114891268240-39741644800E_{p-3})p^4\mod {p^5}.
\end{align*}
\par{\it Proof.} From Corollary 3.5 and Lemma 4.10,
\begin{align*}&\sum_{n=0}^{p-1}\big(1095183237375n^4
+4415024793600n^3+6644229277440n^2
\\&\quad\quad+4436750566144n+1110021214352\big)G_n^2
\\&=16974593p^4\left(G_p^2+65536G_{p-1}^2\right)
-p^2\left(1082146816p^2-1073725440p+397416448\right)G_pG_{p-1}
\\&\e 16974593p^4(12^2+65536\cdot256^{2(p-1)})-p^2
\left(1082146816p^2-1073725440p
+397416448\right)
\\&\q\times\big(12+64(-1)^{\f{p-1}2}p^2E_{p-3}\big)\big( (-1)^{\f{p-1}2}256^{p-1}+3p^2E_{p-3}\big)
\\&\e 1114891268240p^4-12p^2\left(1082146816p^2-1073725440p
+397416448\right)(-1)^{\f{p-1}2}256^{p-1}
\\&\q-39741644800p^4E_{p-3}\mod {p^5}.
\end{align*}
This proves the theorem.
\vskip0.2cm

\par{\bf Conjecture 4.1} Let $p>3$ be a prime and $m,r\in\Bbb Z^+$. Then
\begin{align*}
&S_{mp^r-1}\e (-1)^{\f{p-1}2}32^{mp^{r-1}(p-1)}S_{mp^{r-1}-1}+
\f 14m^2S_mp^{2r}E_{p-3}\mod {p^{2r+1}},
\\&a_{mp^r-1}\e (-1)^{[\f p3]}9^{mp^{r-1}(p-1)}a_{mp^{r-1}-1}+\f 13m^2a_mp^{2r}U_{p-3}\mod {p^{2r+1}},
\\&W_{mp^r-1}\e (-1)^{[\f p3]}27^{mp^{r-1}(p-1)}W_{mp^{r-1}-1}-\f 13m^2W_mp^{2r}U_{p-3}\mod {p^{2r+1}},
\\&G_{mp^r-1}\e (-1)^{\f{p-1}2}256^{mp^{r-1}(p-1)}G_{mp^{r-1}-1}
+\f 14m^2G_mp^{2r}E_{p-3}\mod {p^{2r+1}},
\\&G_{mp^r-1}^{(3)}\e (-1)^{[\f p3]}729^{mp^{r-1}(p-1)}G_{mp^{r-1}-1}^{(3)}
+\f 13m^2G_m^{(3)}p^{2r}U_{p-3}\mod {p^{2r+1}},
\\&G_{mp^r-1}^{(4)}\e (-1)^{[\f p4]}4096^{mp^{r-1}(p-1)}G_{mp^{r-1}-1}^{(4)}
+\f 14m^2G_m^{(4)}p^{2r}s_{p-3}\mod {p^{2r+1}} \ (p\not=5),
\\&G_{mp^r-1}^{(6)}\e (-1)^{\f{p-1}2}186624^{mp^{r-1}(p-1)}G_{mp^{r-1}-1}^{(6)}
+\f 5{108}m^2G_m^{(6)}p^{2r}E_{p-3}\mod {p^{2r+1}},
\\&Q_{mp^r-1}\e (-1)^{[\f p3]}72^{mp^{r-1}(p-1)}Q_{mp^{r-1}-1}
-\f 5{12}m^2Q_mp^{2r}U_{p-3}\mod {p^{2r+1}},
\end{align*}
where $s_n$ is given by $s_0=1$ and
$s_n=1-\sum_{k=0}^{n-1}\b nk2^{2n-1-2k}s_k\ (n\ge 1)$.
\vskip0.2cm
\par{\bf Remark 4.1} The weak version of Conjecture 4.1 is published in [22]. Let $p>3$ be a prime and $m,r\in\Bbb Z^+$. It is worth stating the following challenging conjecture in [22]:
  \begin{align*}
  &S_{mp^r}\e S_{mp^{r-1}}+8m^2S_{m-1}(-1)^{\f{p-1}2\cdot r}p^{2r}E_{p-3}\mod{p^{2r+1}},
  \\&a_{mp^r}\e a_{mp^{r-1}}+3m^2a_{m-1}(-1)^{[\f p3]r}p^{2r}U_{p-3}\mod {p^{2r+1}},
  \\&W_{mp^r}\e W_{mp^{r-1}}-9m^2W_{m-1}(-1)^{[\f p3]r}p^{2r}U_{p-3}\mod {p^{2r+1}},
  \\&G_{mp^r}\e G_{mp^{r-1}}+64m^2G_{m-1}(-1)^{\f{p-1}2\cdot r}p^{2r}E_{p-3}\mod {p^{2r+1}},
  \\&G_{mp^r}^{(3)}\e G_{mp^{r-1}}^{(3)}+243m^2G_{m-1}^{(3)}(-1)^{[\f p3]r}p^{2r}U_{p-3}\mod {p^{2r+1}},
  \\&G_{mp^r}^{(4)}\e G_{mp^{r-1}}^{(4)}+1024m^2G_{m-1}^{(4)}(-1)^{[\f p4]r}p^{2r}s_{p-3}\mod {p^{2r+1}}\ (p\not=5),
  \\&G_{mp^r}^{(6)}\e G_{mp^{r-1}}^{(6)}+8640m^2G_{m-1}^{(6)}(-1)^{\f{p-1}2\cdot r}p^{2r}E_{p-3}\mod {p^{2r+1}},
  \\&Q_{mp^r}\e Q_{mp^{r-1}}-30m^2Q_{m-1}(-1)^{[\f p3]r}p^{2r}U_{p-3}\mod {p^{2r+1}}.
  \end{align*}
In the case $r=1$, the congruence for $S_{mp}-S_m$ modulo $p^3$ was given in [3],
the congruence for $a_{mp}-a_m$ modulo $p^3$ was given in [7],
and the congruence for $Q_{mp}-Q_m$ modulo $p^3$ was given in [14].
\par\q
\par{\bf Lemma 4.11 ([25, Lemma 2.5])}
{\sl For any odd prime $p$ and $b,c\in\Bbb Z$,
$$T_p(b,c)\e b\mod p\q\t{and}\q T_{p-1}(b,c)\e \ls{b^2-4c}p
\mod p.$$}
\par{\bf Theorem 4.11} Let $p$ be an odd prime, $b,c\in\Bbb Z$ and $d=b^2-4c\not\e 0\mod p$. Then for any rational $p$-integer $x$ with $x\not\e 0,-1\mod p$,
\begin{align*}&\sum_{n=0}^{p-1}\Big(\Big(\frac{1-x^2}x
 +\frac{4(x-1)}{x+1}\cdot
 \frac{b^2}{d}\Big)n^2
 +\Big(\frac 2x+\frac{4(x^2-2x-1)}{(x+1)^2}\cdot\frac {b^2}{d}\Big)n
 \\&\quad+ \frac 1x+\frac{(x^2-4x-1)}{(x+1)^2}\cdot
 \frac{b^2}{d}\Big)\frac{T_n(b,c)^2}
 {(dx)^n}
 \\&\e -\f{2b^2(x-1)}{d(x+1)^2}\Ls dpp\mod {p^2}.
 \end{align*}
 \par{\it Proof.} From Corollary 2.1 and Lemma 4.11,
\begin{align*}&\sum_{n=0}^{p-1}\Big(\Big(\frac{1-x^2}x
 +\frac{4(x-1)}{x+1}\cdot
 \frac{b^2}{d}\Big)n^2
 +\Big(\frac 2x+\frac{4(x^2-2x-1)}{(x+1)^2}\cdot\frac {b^2}{d}\Big)n
 \\&\quad+ \frac 1x+\frac{(x^2-4x-1)}{(x+1)^2}\cdot
 \frac{b^2}{d}\Big)\frac{T_n(b,c)^2}
 {(dx)^n}
 \\&=\frac p{(dx)^{p-1}}\Big(p\frac{T_p(b,c)^2}{d}+
 \frac{p}xT_{p-1}(b,c)^2
 -\frac{2b(2p(x+1)+x-1)}{d(x+1)^2}
 T_p(b,c)T_{p-1}(b,c)\Big)
 \\&\e -\f{2b^2(x-1)}{d(x+1)^2}\Ls dpp\mod {p^2}.
 \end{align*}
 \par{\bf Remark 4.2} In [25], Z.W. Sun obtained some congruences for $\sum_{n=0}^{p-1}\f{T_n(b,c)}{m^n}$ modulo
 $p$ or $p^2$, where $p$ is an odd prime.

\end{document}